\documentclass{amsart}
\usepackage{amsfonts}
\usepackage{graphicx}
\usepackage{amscd}

\newtheorem{theorem}{Theorem}[section]
\theoremstyle{plain}

\newtheorem{corollary}[theorem]{Corollary}

\theoremstyle{remark}

\theoremstyle{example}

\numberwithin{equation}{section}

\begin{document}
\title[Minimal graphs]{On the Jacobian of minimal graphs in $\mathbb{R}^{4}$}
\author{Th. Hasanis}
\address{Department of Mathematics, University of Ioannina, 45110, Ioannina,
Greece}
\email{thasanis@uoi.gr}
\author{A. Savas-Halilaj}
\address{Department of Mathematics and Statistics, University of Cyprus, P.O. Box 20537, 1678, Nicosia, 
Cyprus}
\email{halilaj.andreas@ucy.ac.cy}
\author{Th. Vlachos}
\address{Department of Mathematics, University of Ioannina, 45110, Ioannina,
Greece}
\email{tvlachos@uoi.gr}
\thanks{This work was written during the second author's stay at the University of Ioannina  as a research fellow of the ''Foundation for
Education and European Culture''.}

\begin{abstract}
We provide a characterization for complex analytic curves among two-dimensional minimal graphs
in $\mathbb{R}^{4}$ via the Jacobian.
\end{abstract}

\maketitle

\section{Introduction}

The classical theorem of S. Bernstein \cite{B} states that the only entire
minimal graphs in the Euclidean space $\mathbb{R}^{3}$ are planes.
Equivalently, if $u:\mathbb{R}^{2}\rightarrow \mathbb{R}$ is an entire
smooth solution of the differential equation%
\begin{equation*}
{\mbox {div}} \left( \frac{{\mbox {grad}}u}{\sqrt{1+\left\vert {\mbox {grad}}%
u\right\vert ^{2}}}\right) =0,
\end{equation*}%
then $u$ is an affine function. It was conjectured for a long time that the
theorem of Bernstein holds in any dimension. However, for $n=3$,\ its
validity was proved by E. De Giorgi \cite{G}, for $n=4$ by F. Almgren \cite%
{A} and for $n=5,6,7$ by J. Simons \cite{S}. It was a big surprise when E.
Bombieri, E. De Giorgi and E. Giusti \cite{BDG} proved that, for $n\geq 8$,
there are entire solutions of the minimal surface equation other than the
affine ones.

In this paper we study minimal surfaces $M^{2}$ which arise
as graphs over vector valued maps $f:\mathbb{R}^{2}\rightarrow \mathbb{R}
^{2} $, $f=\left( f_{1},f_{2}\right) $, that is%
\begin{equation*}
M^{2}=G_{f}:=\left\{ \left( x,y,f_{1}\left( x,y\right) ,f_{2}\left(
x,y\right) \right) \in \mathbb{R}^{4}:\left( x,y\right) \in \mathbb{R}%
^{2}\right\} .
\end{equation*}%
There are plenty of complete minimal graphs in $\mathbb{R}^{4}$, other than
the planes. More precisely, if $f:\mathbb{C\rightarrow C}$ is any entire
holomorphic or anti-holomorphic function, then the graph $G_{f}$ of $f$ in $%
\mathbb{C}^{2}=\mathbb{R}^{4}$ is a minimal surface and is called a \textit{%
complex analytic curve}. It should be noticed that R. Osserman in \cite{O},
has constructed examples of complete minimal two dimensional graphs in $%
\mathbb{R}^{4}$, which are not complex analytic curves with respect to any
orthogonal complex structure on $\mathbb{R}^{4}$. For instance, the
graph $G_{f}$ over the map $f:\mathbb{R}^{2}\rightarrow \mathbb{R}^{2}$,
given by%
\begin{equation*}
f\left( x,y\right) =\frac{1}{2}\left( e^{x}-3e^{-x}\right) \left( \cos \frac{%
y}{2},-\sin \frac{y}{2}\right) ,\quad \left( x,y\right) \in \mathbb{R}^{2},
\end{equation*}
is such an example. It is worth noticing that the \textit{Jacobian} $%
J_{f}:=\det \left( df\right) $ of $f$  given by $J_{f}=-(e^{2x}-9e^{-2x})/8$ takes every real value.

The problem that we deal with in this article, is to find under which geometric
conditions the minimal\ graph of an entire vector valued map $f:\mathbb{R}%
^{2}\rightarrow \mathbb{R}^{2}$, $f=\left( f_{1},f_{2}\right) $, is a
complex analytic curve. The first result was obtained by S.S. Chern and R.
Osserman \cite{CO}, where they proved that if the differential $df$ of $f$
is bounded, then $G_{f}$ must be a plane. A few years later, L. Simon \cite{Si}
obtained a much more general result by proving that if one of $f_{1}$, $%
f_{2} $ has bounded gradient, then $f$ is affine. Later on, R. Schoen \cite
{Sc} obtained a Bernstein type result by imposing the assumption that $f:%
\mathbb{R}^{2}\rightarrow \mathbb{R}^{2}$ is a diffeomorphism. Moreover, L.
Ni \cite{N} has derived a result of Bernstein type under the assumption that 
$f$ is an \textit{area-preserving map}, that is the Jacobian $J_{f}$  satisfies $J_{f}=1$. This result was
generalized by the authors in a previous paper \cite{HSV}, just by assuming
that $J_{f}$ is bounded.

Another interesting class of minimal surfaces in $\mathbb{R}^{4}$ can be
obtained by considering graphs over vector valued maps of the form $f={\mbox {grad}}u$, where $u:U\subset \mathbb{R}^{2}\rightarrow \mathbb{R}$ is a smooth
function and $U$ an open subset of $\mathbb{R}^{2}$. It can be shown that
the graph over $f={\mbox {grad}}u:U\rightarrow \mathbb{R}^{2}$ is a minimal
surface if and only if the function $u$ satisfies the so called \textit{%
Special Lagrangian Equation}%
\begin{equation*}
\cos \theta \Delta u=\sin \theta \left( \det {\mbox{Hess}} u-1\right) ,
\end{equation*}%
for some real constant$\ \theta $. L. Fu \cite{F}, proved that the entire
solutions of Special Lagrangian Equation are only the harmonic functions or
the quadratic polynomials, which means that the entire minimal graph of $%
{\mbox {grad}}u$ is either a complex analytic curve or a plane.

In this paper, we prove the following result of Bernstein type, from which known results due to R. Schoen \cite
{Sc}, L. Fu \cite{F}, L. Ni \cite{N},  and the authors \cite{HSV} follow as a consequence.

\begin{theorem}
Let $f:\mathbb{R}^{2}\mathbb{\rightarrow R}^{2}$ be an entire smooth vector
valued map such that its graph $G_{f}$ is a minimal surface in $\mathbb{%
R}^{4}$. Assume that $G_{f}$ is not a plane. Then, the graph $G_{f}$ over $f$
is a complex analytic curve if and only if the Jacobian $J_{f}$ of $f$ does
not take every real value. In particular, if $G_{f}$ is a complex analytic curve, then the Jacobian $J_{f}$ takes every real value in $%
\left( 0,+\infty \right) $ or in $\left[ 0,+\infty \right) $ (resp. $%
\left(-\infty, 0 \right) $ or in  $\left(-\infty, 0 \right]$), if $f$ is
holomorphic (resp. anti-holomorphic). 
\end{theorem}

\section{Preliminaries}

At first we set up the notation and recall some basic facts about minimal
surfaces in $\mathbb{R}^{4}$. Let $U$ be an open subset of $\mathbb{R}^{2}$
and $f:U\rightarrow \mathbb{R}^{2}$, $f=\left( f_{1},f_{2}\right) $, a
smooth vector valued map. Then, its graph can be represented by the map $%
X:U\rightarrow \mathbb{R}^{4}$, given by%
\begin{equation*}
X\left( x,y\right) =\left( x,y,f_{1}\left( x,y\right) ,f_{2}\left(
x,y\right) \right) ,\quad \left( x,y\right) \in U.
\end{equation*}%
Denote by $g_{11}$, $g_{12}=g_{21}$ and $g_{22}$ the \textit{coefficients of the
first fundamental form}, which\ are given by%
\begin{equation*}
g_{11}=1+\left\vert f_{x}\right\vert ^{2},\quad g_{12}=g_{21}=\left\langle
f_{x},f_{y}\right\rangle ,\quad g_{22}=1+\left\vert f_{y}\right\vert ^{2},
\end{equation*}%
where $\left\langle \text{ },\text{ }\right\rangle $ stands for the
Euclidean inner product.

Suppose, now, that $\xi $ is a unit vector field  normal to the
surface $G_{f}$. Let $b_{11}\left( \xi \right) $, $b_{12}\left( \xi \right)=b_{21}\left( \xi \right) $
and $b_{22}\left( \xi \right) $ be the \textit{coefficients of the second
fundamental form} with respect to the direction $\xi $, that is%
\begin{equation*}
b_{11}\left( \xi \right) =\left\langle X_{xx},\xi \right\rangle ,\quad
b_{12}\left( \xi \right)=b_{21}\left( \xi \right) =\left\langle X_{xy},\xi \right\rangle ,\quad
b_{22}\left( \xi \right) =\left\langle X_{yy},\xi \right\rangle .
\end{equation*}%
The \textit{mean curvature} $H\left( \xi \right) $, with
respect to the normal direction $\xi $, is defined by the formula%
\begin{equation*}
H\left( \xi \right) :=\frac{%
g_{22}b_{11}\left( \xi \right) -2g_{12}b_{12}\left( \xi \right)
+g_{11}b_{22}\left( \xi \right) }{2\left( g_{11}g_{22}-g_{12}^{2}\right) }.
\end{equation*}%
The surface $G_{f}$ is called \textit{minimal} if it is a critical point of
the area functional. It can be proved that $G_{f}$ is minimal if and only if
the mean curvature, with respect to any unit normal vector
field along $G_{f}$, is zero. A simple computation shows that minimality of $%
G_{f}$ is expressed by the differential equation%
\begin{equation}
\left( 1+\left\vert f_{y}\right\vert ^{2}\right) f_{xx}-2\left\langle
f_{x},f_{y}\right\rangle f_{xy}+\left( 1+\left\vert f_{x}\right\vert
^{2}\right) f_{yy}=0,
\end{equation}%
the so called \textit{minimal surface equation.}

It is a fact of central importance that any surface can be parametrized (at
least locally) by isothermal parameters. This means that there exists a
local diffeomorphism $F\left( u,v\right) =\left( x\left( u,v\right) ,y\left(
u,v\right) \right) $ such that the coefficients of the first fundamental
form of the new parametrization%
\begin{equation*}
X\left( u,v\right) =\left( x\left( u,v\right) ,y\left( u,v\right)
,f_{1}\left( x\left( u,v\right) ,y\left( u,v\right) \right) ,f_{2}\left(
x\left( u,v\right) ,y\left( u,v\right) \right) \right)
\end{equation*}%
become%
\begin{equation*}
g_{11}\left( u,v\right) =g_{22}\left( u,v\right) =E\left( u,v\right) \text{
\ and \ }g_{12}\left( u,v\right) =g_{21}\left( u,v\right)=0.
\end{equation*}
Set
\begin{equation*}
\varphi \left( u,v\right): =f_{1}\left( x\left( u,v\right) ,y\left(
u,v\right) \right) \text{ \ and \ }\psi \left( u,v\right): =f_{2}\left(
x\left( u,v\right) ,y\left( u,v\right) \right) .
\end{equation*}
Then, in isothermal parameters, the minimal surface equation is equivalent to
\begin{equation*}
x _{uu}+x _{vv}=0=y _{uu}+y _{vv},
\end{equation*}
\begin{equation*}
\varphi _{uu}+\varphi _{vv}=0=\psi _{uu}+\psi _{vv}.
\end{equation*}
Furthermore, the complex valued functions $\phi _{k}:U\mathbb{\subseteq C\rightarrow C}$, $k=1,2,3,4,$ given by
\begin{equation*}
\left\{ 
\begin{array}{ll}
\phi _{1}=x_{u}-ix_{v}, & \phi _{2}=y_{u}-iy_{v}, \\ 
\phi _{3}=\varphi _{u}-i\varphi _{v}, & \phi _{4}=\psi _{u}-i\psi _{v}%
\end{array}%
\right.
\end{equation*}%
are holomorphic and satisfy%
\begin{equation*}
\phi _{1}^{2}+\phi _{2}^{2}+\phi _{3}^{2}+\phi _{4}^{2}=0.
\end{equation*}

Assume, now, that $f:\mathbb{R}^{2}\rightarrow \mathbb{R}^{2}$ is an entire
solution of the minimal surface equation. In this case, R. Osserman \cite[%
Theorem 5.1]{O} proved the following result which will be the main tool for
the proof of our results.

\begin{theorem}
Let $f:\mathbb{R}^{2}\rightarrow \mathbb{R}^{2}$ be an entire solution of
the minimal surface equation. Then there exist real constants $a$, $b$, with 
$b>0$, and a non-singular transformation%
\begin{equation*}
x=u,\quad y=au+bv,
\end{equation*}%
such that $\left( u,v\right) $ are global isothermal parameters for the
surface $G_{f}$.
\end{theorem}

\section{Proof of the theorem}

Let $f=\left( f_{1},f_{2}\right) $ be an entire solution of the minimal
surface equation%
\begin{equation*}
\left( 1+\left\vert f_{y}\right\vert ^{2}\right) f_{xx}-2\left\langle
f_{x},f_{y}\right\rangle f_{xy}+\left( 1+\left\vert f_{x}\right\vert
^{2}\right) f_{yy}=0.
\end{equation*}%
Then its graph,%
\begin{equation*}
G_{f}=\left\{ \left( x,y,f_{1}\left( x,y\right) ,f_{2}\left( x,y\right)
\right) \in \mathbb{R}^{4}:\left( x,y\right) \in \mathbb{R}^{2}\right\} ,
\end{equation*}%
is a minimal surface. According to Theorem 2.1, we can
introduce global isothermal parameters, via a non-singular transformation%
\begin{equation*}
x=u,\text{ \ }y=au+bv,
\end{equation*}%
where $a$, $b$ are real constants with $b>0$. Now, the minimal surface $%
G_{f} $ is parametrized by the map%
\begin{equation*}
X\left( u,v\right) =\left( u,au+bv,\varphi \left( u,v\right) ,\psi \left(
u,v\right) \right),
\end{equation*}%
where%
\begin{equation*}
\varphi \left( u,v\right): =f_{1}\left( u,au+bv\right) \text{ \ and \ }\psi
\left( u,v\right): =f_{2}\left( u,au+bv\right) .
\end{equation*}%
Set $\Phi =\left( \varphi ,\psi \right) $. Because of the relation%
\begin{equation*}
\frac{\partial \left( \varphi ,\psi \right) }{\partial \left( u,v\right) }=%
\frac{\partial \left( f_{1},f_{2}\right) }{\partial \left( x,y\right) }\frac{%
\partial \left( x,y\right) }{\partial \left( u,v\right) }
\end{equation*}%
for the Jacobians, we obtain%
\begin{equation*}
J_{\Phi }=bJ_{f}.
\end{equation*}%
Since $\left( u,v\right) $ are isothermal parameters and taking into account
that $G_{f}$ is minimal, it follows that the functions $\varphi $ and $\psi $
are harmonic, that is%
\begin{equation*}
\varphi _{uu}+\varphi _{vv}=0=\psi _{uu}+\psi _{vv}.
\end{equation*}%
Then, the complex valued functions $\phi _{k}:\mathbb{C\rightarrow C}$, $%
k=1,2,3,4,$ given by%
\begin{equation}
\left\{ 
\begin{array}{ll}
\phi _{1}=1, & \phi _{2}=a-ib, \\ 
\phi _{3}=\varphi _{u}-i\varphi _{v}, & \phi _{4}=\psi _{u}-i\psi _{v}%
\end{array}%
\right.
\end{equation}%
are holomorphic and satisfy%
\begin{equation}
\phi _{1}^{2}+\phi _{2}^{2}+\phi _{3}^{2}+\phi _{4}^{2}=0.
\end{equation}

\begin{proof}[Proof of Theorem 1.1]
Assume that the graph $G_{f}$ of $f(x,y)=(f_{1}(x,y),f_{2}(x,y))$, $(x,y) \in \mathbb{R}^{2}$, is a minimal surface which is not a plane.  

Suppose now that $J_{f}$ does not take every real value. We
shall prove that $G_{f}$ is a complex analytic curve. Arguing indirectly,
we assume that $G_{f}$ is not a complex analytic curve. Then, 
equation 
$\left( 3.2\right)$ can be written equivalently in the form%
\begin{equation}
\left( \phi _{3}-i\phi _{4}\right) \left( \phi _{3}+i\phi _{4}\right) =-d,
\end{equation}%
where%
\begin{equation*}
d=1+\left( a-ib\right) ^{2}.
\end{equation*}%
We claim that $d\neq 0$. Assume to the contrary that $d=0$. Then, $a=0$, $%
b=1 $ and consequently $\left( x,y\right)$ are isothermal parameters. Furthermore, from $%
\left( 3.3\right) $ we obtain $\varphi _{3}=\pm i\varphi _{4}$, or
equivalently,%
\begin{equation}
\frac{\partial f_{1}}{\partial x}-i\frac{\partial f_{1}}{\partial y}=\pm
i\left( \frac{\partial f_{2}}{\partial x}-i\frac{\partial f_{2}}{\partial y}%
\right) .
\end{equation}%
From $\left( 3.4\right) $ we readily deduce that $f=f_{1}+if_{2}$ is
holomorphic or anti-holomorphic, which is\ a contradiction. Hence $d\neq 0$.

By virtue of $\left( 3.3\right) $, we see that $\varphi _{3}-i\varphi
_{4} $, $\varphi _{3}+i\varphi _{4}$ are entire nowhere vanishing
holomorphic functions. Define the complex valued function $h:\mathbb{%
C\rightarrow C}$, by%
\begin{equation}
h=\phi _{3}-i\phi _{4}.
\end{equation}%
We point out that $h$ is holomorphic, non-constant and nowhere vanishing.
Combining $\left( 3.3\right) $ with $\left( 3.5\right) $, we get 
\begin{equation}
\phi _{3}=\frac{1}{2}\left( h-\frac{d}{h}\right) \text{ \ and \ }\phi _{4}=%
\frac{i}{2}\left( h+\frac{d}{h}\right) .
\end{equation}%
Bearing in mind $\left( 3.1\right) $, it follows that the imaginary part of $\phi _{3}\overline{\phi }_{4}$ is given by
\begin{equation*}
{\mbox{Im}}(\phi _{3}\overline{\phi }_{4})=\varphi _{u}\psi _{v}-\varphi
_{v}\psi _{u}=J_{\Phi }.
\end{equation*}%
On the other hand, from $\left( 3.6\right) $ we get%
\begin{equation*}
{\mbox{Im}}(\phi _{3}\overline{\phi }_{4})=\frac{1}{4}\left( -\left\vert
h\right\vert ^{2}+\frac{\left\vert d\right\vert ^{2}}{\left\vert
h\right\vert ^{2}}\right) .
\end{equation*}%
Thus, taking into account the relation $J_{\Phi }=bJ_{f}$, we have
\begin{equation*}
J_{f}=\frac{1}{4b}\left( -\left\vert h\right\vert ^{2}+\frac{\left\vert
d\right\vert ^{2}}{\left\vert h\right\vert ^{2}}\right) .
\end{equation*}%
Since $h$ is an entire and non-constant holomorphic function, by Picard's
Theorem, there are sequences $\left\{ z_{n}\right\} _{n\in \mathbb{N}}$ and $%
\left\{ w_{n}\right\} _{n\in \mathbb{N}}$ of complex numbers\ such that $%
\left\vert h\left( z_{n}\right) \right\vert \rightarrow \infty $ and $%
\left\vert h\left( w_{n}\right) \right\vert \rightarrow 0$. Consequently, $J_{f}\left(
z_{n}\right) \rightarrow-\infty $ and $J_{f}\left( w_{n}\right) \rightarrow\infty $. Thus, $%
J_{\Phi }\left( \mathbb{R}^{2}\right) =\mathbb{R}$, which contradicts to our
assumptions. Therefore, $G_{f}$ must be a complex analytic curve.

Conversely, assume that $G_{f}$ is a complex analytic curve which is not a plane. Then, the complex
valued function $f=f_{1}+if_{2}$ is holomorphic or anti-holomorphic. We
introduce the complex variable $z=x+iy$. An easy computation shows that%
\begin{equation*}
J_{f}=\left\vert f_{z}\right\vert ^{2}-\left\vert f_{\overline{z}%
}\right\vert ^{2}.
\end{equation*}%

Hence, $J_{f}\geq 0$ if $f$ is holomorphic and $J_{f}\leq 0$ if $f$ is
anti-holomorphic. In either case, $J_{f}$ does not take every real value.

Asssume at first that $f$ is holomorphic. Then, 
\begin{equation*}
J_{f}=\left\vert f_{z}\right\vert ^{2}.
\end{equation*}%
By minimality, we obtain that $f_{z}$ is an entire holomorphic function.
Moreover $f_{z}$ cannot be constant, since otherwise $f$ is affine and $%
G_{f} $ a plane. Consequently, appealing to Picard's Theorem, the range of $f_{z}$ is the whole complex plane $\mathbb{C}$, or the plane minus a single point. Thus, the range of $%
J_{f} $ must be $\left( 0,+\infty \right) $ or $\left[ 0,+\infty \right) $.

Suppose, now, that $f$ is anti-holomorphic. Then,%
\begin{equation*}
J_{f}=-\left\vert f_{\overline{z}}\right\vert ^{2}
\end{equation*}%
and $f_{\overline{z}}$ is an entire anti-holomorphic function. Arguing as
above, we deduce that range of $J_{f}$ must be $\left( -\infty ,0\right) $
or $\left( -\infty ,0\right] .$
\end{proof}

\section{Applications}

In this last section, we reobtain some known Bernstein type theorems for entire vector valued maps $f:%
\mathbb{R}^{2}\rightarrow \mathbb{R}^{2}$ using the method developed here. The following result was first
obtained by L. Simon \cite{Si}.

\begin{corollary}
Let $f:\mathbb{R}^{2}\rightarrow \mathbb{R}^{2}$, $f=\left(
f_{1},f_{2}\right) $ be an entire solution of the minimal surface equation,
such that one of $f_{1}$, $f_{2}$ has bounded gradient. Then, $f$ is an
affine map.
\end{corollary}

\begin{proof}
Without loss of generality, let us assume that $f_{1}$ has bounded gradient. By
virtue of Theorem $2.1$, we can introduce global isothermal
parameters $\left( u,v\right) $, via the non-singular transformation%
\begin{equation*}
x=u,\text{ \ }y=au+bv,
\end{equation*}%
where $a$, $b$ are real constants with $b>0$. Then, the graph $G_{f}$ is
parametrized via the map%
\begin{equation*}
X\left( u,v\right) =\left( u,au+bv,f_{1}\left( u,au+bv\right) ,f_{2}\left(
u,au+bv\right) \right) .
\end{equation*}%
Consider the real valued functions $\varphi $, $\psi :\mathbb{R}%
^{2}\rightarrow \mathbb{R}$, given by%
\begin{equation*}
\varphi \left( u,v\right) :=f_{1}\left( u,au+bv\right) ,\text{ \ }\psi
\left( u,v\right) :=f_{2}\left( u,au+bv\right) .
\end{equation*}%
Since $\left( u,v\right) $ are isothermal parameters, $\varphi $, $\psi $ are
harmonic functions. Moreover, the complex valued functions $\phi _{k}:%
\mathbb{C\rightarrow C}$, $k=1,2,3,4,$ given by%
\begin{equation*}
\left\{ 
\begin{array}{ll}
\phi _{1}=1, & \phi _{2}=a-ib, \\ 
\phi _{3}=\varphi _{u}-i\varphi _{v}, & \phi _{4}=\varphi _{u}-i\psi _{v},%
\end{array}%
\right.
\end{equation*}%
are holomorphic and satisfy%
\begin{equation}
\phi _{1}^{2}+\phi _{2}^{2}+\phi _{3}^{2}+\phi _{4}^{2}=0.
\end{equation}%
A simple calculation shows that%
\begin{equation*}
\frac{\partial \varphi }{\partial u}=\frac{\partial f_{1}}{\partial x}+a%
\frac{\partial f_{1}}{\partial y},\text{ \ }\frac{\partial \varphi }{%
\partial v}=b\frac{\partial f_{1}}{\partial y}.
\end{equation*}%
Since $f_{1}$ has bounded gradient, it follows that $\varphi $ has also
bounded gradient. Because $\varphi _{u}$ and $\varphi _{v}$ are harmonic
functions, by Liouville's Theorem it follows that they must be constants.
Hence $\phi _{3}$ is a constant complex valued function. By virtue of  $\left(
4.1\right) $, we obtain that $\phi _{4}$ must also be constant. Hence $G_{f}$
is a plane and $f$ is an affine function.
\end{proof}

The following result due to R. Schoen \cite{Sc} is a consequence
of Theorem 1.1.

\begin{corollary}
Let $f:\mathbb{R}^{2}\rightarrow \mathbb{R}^{2}$ be an entire solution of
the minimal surface equation. If $f$ is a diffeomorphism, then $f$ is an
affine map.
\end{corollary}

\begin{proof}
Since $f$ is a diffeomorphism, it follows that $J_{f}>0$ if $f$ is
orientation preserving, or $J_{f}<0$ if $f$ is orientation reversing. Then,
by Theorem 1.1, $f$ must be holomorphic or anti-holomorphic. Thus, $f$ is an
entire conformal diffeomorphism. By a classical theorem of Complex Analysis 
\cite[p. 388]{P}, $f$ must be an affine map.
\end{proof}

Finally, we provide an alternative proof of the result due to L. Fu \cite{F},  based on Theorem 1.1.

\begin{corollary}
Let $u:\mathbb{R}^{2}\rightarrow \mathbb{R}$ be an entire solution of the
Special Lagrangian Equation%
\begin{equation*}
\cos \theta \Delta u=\sin \theta \left( \det {\rm {Hess}}u-1\right) ,
\end{equation*}%
where $\theta $ is a real constant. Then, either $u$ is a harmonic function
or $u$ is a quadratic polynomial.
\end{corollary}

\begin{proof}
Consider the entire vector valued map $f={\mbox{grad}}u:\mathbb{R}^{2}\to\mathbb{R}^{2}$. Since $u$ satisfies the Special Lagrangian Equation, the graph of $f$
is a minimal surface in $\mathbb{R}^{4}$. Note that in this case%
\begin{equation*}
J_{f}=\det {\mbox{Hess}}u=u_{xx}u_{yy}-u_{xy}^{2}.
\end{equation*}%
Suppose at first that there is a point $\left( x_{0},y_{0}\right) \in 
\mathbb{R}^{2}$ such that $J_{f}\left( x_{0},y_{0}\right) =1$. Then, at this
point the Laplacian of $u$  satisfies
\begin{equation*}
\Delta u\left( x_{0},y_{0}\right) =u_{xx}\left( x_{0},y_{0}\right)
+u_{yy}\left( x_{0},y_{0}\right) \neq 0.
\end{equation*}%
Consequently, $\theta =\frac{\pi }{2}$ and $J_{f}\equiv 1$. By Theorem 1.1, $f$ must be affine and thus $u$ a quadratic polynomial.

Assume now that $J_{f}\neq 1$. Then, either $J_{f}>1$ or $J_{f}<1$. If $J_{f}>1$,
then, according to Theorem 1.1, $f$ must be an affine map and so $u$ is a quadratic polynomial. If $%
J_{f}<1$, then, by virtue of Theorem 1.1, we deduce that  $J_{f}\leq 0$ and $f$ is an anti-holomorphic function.
Appealing to Cauchy-Riemann equations, it follows that $u$ is a harmonic
function.
\end{proof}

\end{document}